\documentclass{tac}

\usepackage{amsmath, amssymb}
\usepackage{enumitem}
\usepackage{bbold}
\usepackage{xy}
\input diagxy

\usepackage[colorlinks=true]{hyperref}
\hypersetup{allcolors=[rgb]{0.1,0.1,0.4}}

\numberwithin{equation}{section}
 
\author{Hongliang Lai and Qingzhu Luo}

\thanks{We would like to thank the support of the National Natural Science Foundation of China (No. 12171342), and Sichuan National Applied Mathematics Center- Chengdu University of Information Technology, Institute of Applied Mathematics for Intelligent Systems (No. 2025ZX003)
}

\address{School of Mathematics, Sichuan University, Chengdu 610064, China. \\College of Applied Mathematics, Chengdu University of Information Technology, Chengdu 610225, China.}

\title{On
the Cartesian closedness of [0,1]-{\sf Cat} and some of its subcategories}

\copyrightyear{2020}


\eaddress{hllai@scu.edu.cn\CR qingzhuluo@cuit.edu.cn}

\newtheorem{thm}{Theorem}

\newtheorem{lem}{Lemma}
\newtheorem{prop}{Proposition}

\newtheoremrm{rem}{Remark}
\newtheoremrm{defn}{Definition}
\newtheoremrm{exmp}{Example}
\newtheoremrm{ques}{Question}

\mathrmdef{Hom}
\mathbfdef{Set}

\begin{document}
\bibliographystyle{abbrv}
\newcommand{\ra}{\rightarrow}
\newcommand{\lra}{\longrightarrow}
\newcommand{\id}{\mathrm{id}}
\newcommand{\bfA}{\mathbf{A}}
\newcommand{\bfB}{\mathbf{B}}
\newcommand{\bfC}{\mathbf{C}}
\newcommand{\sQ}{\mathsf{Q}}
\newcommand{\sM}{\mathsf{M}}
\newcommand{\sK}{\mathsf{K}}
\newcommand{\sN}{\mathsf{N}}
\newcommand{\qcat}{{\sQ\text{-}\mathsf{Cat}}}
\newcommand{\mcat}{\mathsf{M}\text{-}\mathsf{Cat}}
\newcommand{\kcat}{\mathsf{K}\text{-}\mathsf{Cat}}
\newcommand{\sqcat}{{\sQ\text{-}\mathsf{SymCat}}}
\newcommand{\skcat}{{\sK\text{-}\mathsf{SymCat}}}
\newcommand{\qref}{\sQ\text{-}\mathsf{Ref}}
\newcommand{\YCom}{\RCat_{\sf yon}}
\newcommand{\yock}{\mathsf{YonCom(\mathsf{K})}}
\newcommand{\yocm}{\mathsf{YonCom(\mathsf{M})}}
\newcommand{\CCom}{\RCat_{\sf cau}}
\newcommand{\SCom}{\RCat_{\sf smy}}
\newcommand{\ev}{\mathrm{ev}}
\newcommand{\vep}{\varepsilon}

\newcommand{\ordlu}{{ [0,1]_{\L}\text{-}\mathsf{Ord}}}
\newcommand{\ordp}{{ [0,1]_{P}\text{-}\mathsf{Ord}}}
\newcommand{\ord}{{[0,1]\text{-}\mathsf{Ord}}}
\newcommand{\RCat}{[0,1]\text{-}\mathsf{Cat}}
\newcommand{\RfCat}{\RCat_{\sf fin}}
\newcommand{\SRCat}{[0,1]\text{-}\mathsf{SymCat}}
\newcommand{\sqref}{{\sQ\text{-}\mathsf{SymRef}}}
\newcommand{\QRel}{\sQ\text{-}\mathsf{Rel}}
\newcommand{\bbt}{\mathbb{2}}
\newcommand{\Ra}{\Rightarrow}
\newcommand{\sft}{\mathbb{2}}
\newcommand{\Pord}{\mathsf{Ord}}
\newcommand{\bv}{\bigvee}
\newcommand{\bw}{\bigwedge}
\newcommand{\Lra}{\Longrightarrow}

\newcommand{\Idm}{\mathrm{Idm}_{\&}}
\newcommand{\Lu}{\text{\L}}
\newcommand{\lam}{\lambda}
\newcommand{\de}{\delta}
\newcommand{\ga}{\gamma} 
\newcommand{\al}{\alpha}
\newcommand{\be}{\beta}
\newcommand{\with}{\&}
\newcommand{\xlam}{\{x_\lambda\}_{\lambda\in D}}
\newcommand{\flamx}{\{f_{\lambda}(x)\}_{\lambda\in D}}
\newcommand{\fxlam}{\{f(x_{\lambda})\}_{\lambda\in D}}
\newcommand{\flam}{\{f_\lambda\}_{\lambda\in D}}
\newcommand{\ylam}{\{y_\lambda\}_{\lambda\in D}}
\newcommand{\zlam}{\{z_\lambda\}_{\lambda\in D}}
\newcommand{\dt}{d_{\Pi}}





\maketitle

\begin{abstract}
We describe all left continuous triangular norms for which the category $\RCat$ of real-enriched categories and functors is cartesian closed. We
furthermore show that the cartesian closedness of $\RCat$ is equivalent to the cartesian closedness of either (and thus all) of the following subcategories: the full subcategory of Cauchy complete $[0,1]$-categories; the subcategory of Yoneda complete $[0,1]$-categories and Yoneda continuous
$[0,1]$-functors; the full subcategory of Smyth complete $[0,1]$-categories;
and the full subcategory of finite $[0,1]$-categories.
\end{abstract}

\section{Introduction}
In 1973 \cite{Lawvere1973},  it is pointed out by Lawvere that categories enriched over a monoidal closed
category can be viewed as ``ordered sets'' whose truth-values are taken in that closed category. This  viewpoint has led to the  quantitative domain theory, of which the core
objects are categories enriched over a quantale, see e.g.,  \cite{Bonsangue1998,Flagg1997,Hofmann2011,Hofmann2012,Wagner1997}.

Compared with the category $\mathsf{Ord}$ consisting of all ordered sets, failure to be cartesian closed in general is one of the main defects of the category $\qcat$ consisting of categories enriched over a quantale $\sQ$. For instance, if the quantale is $\sQ=([0,1],\&,1)$, that is, the binary operation $\&$ is a left continuous triangular norm on $[0,1]$, then categories enriched over $\sQ$ are called real-enriched categories \cite{Zhang2024}, which are of particular interests in quantitative domain theory. However,  if the left continuous triangular norm $\&$ is also continuous, then all real-enriched categories form a   cartesian closed category if and only if $\&=\wedge$ \cite{Lai2017}. 

The purpose of this paper is to show that dropping the continuity condition of the triangular norm $\&$,  then there are non-trivial left continuous triangular norms different from the operation $\wedge$ such that all real-enriched categories constitute a cartesian closed category. In fact, due to the well-ordering and density of the unit interval, we describe explicitly all such left continuous triangular norms. Moreover, we show that, all real-enriched categories form a  cartesian closed category if and only if all Cauchy complete real-enriched categories  form a cartesian closed full subcategory, if and only if all Yoneda complete real-enriched categories and Yoneda continuous functors form a cartesian closed subcategory, if and only if all Smyth complete real-enriched categories also form a cartesian closed full subcategory.

 The content is arranged as follows: Section 2,   recall some basis about quantale-enriched categories   and the cartesian closedness of the category $\qcat$; Section 3, characterize the left continuous triangular norms such that $\RCat$ is cartesian closed; Section 4 and 5, we show that certain subcategories of $\RCat$ are also cartesian closed, including whose objects are Cauchy complete, Yoneda complete and Smyth complete real-enriched categories.

\section{Quantale-enriched categories}
A \emph{quantale} $\sQ=(Q,\&,k)$ is a complete lattice $Q$ equipped with a monoidal structure  whose binary operation $\&$ preserves all suprema in each place, that is,  
  $$x \& (\bv_{i \in I}{y_i})=\bv_{i \in I}(x \& y_i), \quad
  (\bv_{i\in I}y_i)\& x=\bv_{i\in I}y_i\& x$$ for all $x,y_i\in Q$. 
  The $\&$-neutral element is denoted by $k$. Except in the trivial case, we always have that $k$ is not the bottom element. The quantale $\sQ=(Q,\&,k)$ is \emph{integral} if $k=\top$, where $\top$ is the top element of $Q.$


Let $\sQ=(Q,\&,k)$ be a quantale. A \emph{$\sQ$-category} $(X,r)$ consists of a set (its ``objects'') $X$ and a $\sQ$-relation
 (its ``hom-functor'') $r:X\times X\lra Q$, such that 
 $$k\leq r(x,x),\quad r(y,z)\&r(x,y)\leq r(x,z)$$
 for all $x$, $y$ and $z$ in $X$. 

 A map $f:(X,r)\lra (Y,s)$ between $\sQ$-categories is a \emph{$\sQ$-functor} if it satisfies 
$$r(x,y)\leq s(f(x),f(y))$$
for all $x$ and $y$ in $X$.

All the $\sQ$-categories and $\sQ$-functors constitute a category \[\qcat.\] 



Let $(X,r)$ and $(Y,s)$ be two $\sQ$-categories. The \emph{product} of $(X,r)$ and $(Y,s)$ has the underlying set $X\times Y$ and the $\sQ$-relation $r\times s$ given by \[r\times s((x_1,y_1),(x_2,y_2))=r(x_1,x_2)\wedge s(y_1,y_2)\] for $(x_1,y_1),(x_2,y_2)\in X\times{Y}$. The terminal object $\mathbb{T}$ in $\qcat$ consists of a singleton set $\{\star\}$ and a $\sQ$-relation $e$ with $e(\star,\star)=\top.$

Recall from \cite{Adamek1990} that a category $\mathbf{A}$ is  \emph{cartesian closed} provided that it has finite products and for
each object $A$ the functor $A\times(-): \mathbf{A}\lra\mathbf{A}$ has a right adjoint $(-)^A:
\mathbf{A}\lra\mathbf{A}$. The object $B^A$ is called the power of $A$ and $B$. The category $\qcat$ is cartesian closed if and only if each morphism $t:(X,r)\lra\mathbb{T}$ is \emph{exponentiable} in the sense of \cite{Clementino2006, Clementino2003b}. In this case, we also say that the $\sQ$-category $(X,r)$ is exponentiable for short.   

Notice that the category $\qcat$ is concrete over the category $\mathsf{Set}$, and its terminal object is discrete. If $\qcat$ is cartesian closed, then it has function spaces (see Proposition 27.18 in \cite{Adamek1990}). That is, for $\sQ$-categories $(X,r)$, $(Y,s)$, we can choose  the power object $(Y,s)^{(X,r)}$ with the underlying set consisting of all $\sQ$-functors from $(X,r)$ to $(Y,s)$, denoted by $[(X,r),(Y,s)]$, and the evaluation morphism $$\mathrm{ev}:X\times[(X,r),(Y,s)]\lra Y, \quad\mathrm{ev}(x,f)=f(x).$$

In \cite{Clementino2009b,Stubbe2025}, necessary and sufficient conditions are established for 
quantaloid-enriched categories to be exponentiable, resp. the category of quantaloid-enriched categories 
and functors to be cartesian closed. Below we formulate these conditions in the simpler case where we 
enrich over a quantale (viewed as a quantaloid with a single object).
\begin{prop}{\rm (\cite{Clementino2009b})}\label{exp obj}
A $\sQ$-category $(X,r)$ is exponentiable  if and only if the following two conditions hold:
\begin{enumerate}
\item[{\rm (1)}] for all $\{q_i\mid i\in I\}\subseteq Q$ and all $x,y\in X$, $(\bv_{i\in I}q_i)\wedge r(x,y)=\bv_{i\in I}(q_i\wedge r(x,y))$,
\item[{\rm (2)}] for all $p,q\in Q$ and all $x,z\in X$, 
$(p\&q)\wedge r(x,z)=\bv_{y\in X}(p\wedge r(y,z))\&(q\wedge r(x,y))$.
\end{enumerate}
\end{prop}

\begin{prop}{\rm (\cite{Stubbe2025})}\label{ccc con} The category $\qcat$ is cartesian closed if and only if the following two conditions hold:
\begin{itemize}
\item[{\rm (1)}] the complete lattice $Q$ is a frame,
\item[{\rm (2)}] for all $p,q,u\in Q$,
$(p\&q)\wedge u=((p\wedge u)\&(q\wedge k))\vee((p\wedge k)\&(q\wedge u))$.
\end{itemize}
\end{prop}






\begin{prop} \label{int con} Let $\sQ$ be a quantale with the underlying complete lattice $Q=[0,1]$. If $\qcat$ is cartesian closed, then $k=1$, that is, the quantale is integral, and $p\&p\&p=p\&p$ for all $p\in[0,1]$. 
\end{prop}

\begin{proof} 
Firstly,  for all $p<k$,   by the condition (2) in Proposition \ref{ccc con}, we have that 
$$
(p\&1)\wedge k=((p\wedge k)\&(1\wedge k))\vee((p\wedge k)\&(1\wedge k))=(p\&k)=p.
$$
Thus,  $p\&1=p$  since $p<k$.
Consequently, we see that
$$1=k\&1=(\sup_{p<k}p)\&1=\sup_{p<k}(p\&1)=k.$$
Therefore, the top element $1=k$, hence the quantale is integral. 

Secondly, by the condition (2) in Proposition \ref{ccc con} again, for each $p\in[0,1]$, it holds that
\begin{align*}p\&p &=(p\&p)\wedge(p\&p)\\
                   &=((p\wedge(p\&p))\&(p\wedge k))\vee((p\wedge k)\&(p\wedge(p\&p)))\\
                   &=((p\&p)\&p)\vee(p\&(p\&p))\\
                   &=(p\&p)\&p.
\end{align*}
\end{proof}

\begin{rem}
    In fact, it is shown in \cite{Stubbe2025} that $p\&p\&p=p\&p$ for all $p\in Q$ holds in any integral quantale $\sQ=(Q,\&,k)$ satisfying the condition  (2) in Proposition \ref{ccc con}.
\end{rem}

 \begin{exmp} Let $Q=\{0,k,1\}$ with $0<k<1$. Define a binary operation $\&$ on $Q$ by 
 $$p\& q=\begin{cases} 0& p=0\text{ or } q=0,\\
                       p\vee q & \text{ otherwise}.
        \end{cases}
 $$
 Then $\sQ=(Q,\&,k)$ is a commutative quantale with $k$ being the unit. Clearly,  it satisfies the conditions in Proposition \ref{ccc con}. Thus, the category $\qcat$ is cartesian closed but $\sQ$ is not integral. 
\end{exmp}
 
Given $\sQ$-categories $(X,r)$ and $(Y,s)$, we always choose $$[(X,r),(Y,s)]=\qcat((X,r),(Y,s))$$
as the underlying set of the power object $(Y,s)^{(X,r)}$. 
   When $\sQ$ is integral, the $\sQ$-categorical hom-functor on $[(X,r),(Y,s)]$ is given in the below (see \cite{Clementino2006,Clementino2009}):
  \begin{equation}\label{exp hom}
  d(f,g)=\bv\{q\in Q\mid q\wedge r(x,y)\leq s(f(x),g(y))\text{ for all } x,y\in X\}.
  \end{equation} 
The hom-functor $d$ on function spaces will appear frequently in the rest of this paper.

\section{Real-enriched categories}

From now on, we always consider  a commutative and  integral   quantale with the underlying complete lattice $[0,1]$. In this case, the binary operation $\&$ on $[0,1]$ is called a left continuous triangular norm (t-norm for short) \cite{Klement2000}.

A $\sQ$-category is  called a \emph{real-enriched category} (or \emph{$[0,1]$-category}), a $\sQ$-functor is called a \emph{$[0,1]$-functor} \cite{Zhang2024}.
All real-enriched categories and $[0,1]$-functors constitute a category \[\RCat.\]



Clearly, the complete lattice $[0,1]$ is a frame, thus, it satisfies the first condition in Proposition \ref{ccc con}. But the second condition need not hold for a left continuous t-norm $\&$ on $[0,1]$ in general. 
In the below, we give equivalent formulations of condition (2) in Proposition \ref{ccc con}, specifically for left continuous t-norms, which will be required later on. 
\begin{prop}\label{E1}
Let $\&$ be a left continuous t-norm on $[0,1]$. The following statements are equivalent:
    \begin{enumerate}
    \item[{\rm (1)}] For all $p$, $q$ and $u$ in $[0,1]$,
    \begin{equation}\label{C1}
    (p\&q)\wedge u=((p\wedge u)\&q)\vee(p\& (q\wedge u));
    \end{equation}
    \item[{\rm (2)}] For all $p,u\in[0,1]$,
    \begin{equation}\label{C2}
    u\leq p\& p\Lra u\&p=u;
    \end{equation}
    \item[{\rm (3)}] There is a (countable) family of pairwise disjoint closed intervals $\{[a_i,b_i]\subseteq[0,1)\mid i\in I\}$ such that for all $p,q\in[0,1],$
\begin{equation}\label{C3}
p\& q=\begin{cases} a_i & p,q\in[a_i,b_i]\text{ for some } i\in I,\\
        p\wedge q& otherwise.
\end{cases}
\end{equation}
    \end{enumerate}
\end{prop}

\begin{proof}
(1)$\Lra$(2): Given $p$ and $u$ in $[0,1]$, if   $u\leq p\& p$, then it holds that 
\[u=(p\& p)\wedge u=((p\wedge u)\& p)\vee (p\&(p\wedge u))=(u\&p)\vee(p\&u)=p\&u.\]

(2)$\Lra$(3): 
Collect all idempotent elements 
 of $\&$ in $[0,1]$, we obtain  a subset $Idm$ of $[0,1]$, which is closed under all suprema. For each $a\in Idm$, let $$\hat{a}=\sup\{x\in[0,1]\mid x\& x=a\},$$
 then $\hat{a}$ is the largest element in $[0,1]$ such that $\hat{a}\&\hat{a}=a$ since $\&$ is left continuous. Collect all the idempotent elements $a$ with $a<\hat{a}$, then one obtains a subset $S\subseteq Idm$.
   
For all $a,b\in S$ with $a<b$, one can see that $\hat{a}<b$. Otherwise, if $b\leq\hat{a}$, then $b=b\&b\leq\hat{a}\&\hat{a}=a$, that is a contradiction.  Thus, it follows that $[a,\hat{a}]\cap[b,\hat{b}]=\emptyset$.  That is,  the family  $\{[a,\hat{a}]\mid {a\in S}\}$  is pairwise disjoint, therefore necessarily 
countable, and each interval $[a,\hat{a}]\subseteq[0,1)$. 
 
Given $p$ and $q$ in $[0,1]$ with $p\leq q$, we calculate  the value of $p\&q$ in two cases:
\begin{itemize}
\item[(i)] There is some $a\in S$ such that $a\leq p\leq q\leq\hat{a}$. Clearly, it holds that
$$a=a\&a\leq p\&q\leq\hat{a}\&\hat{a}=a.$$
Thus, $a=p\& q$.

\item[(ii)] For each $a\in S$, either $p\not\in[a,\hat{a}]$ or $q\not\in[a,\hat{a}]$. In fact, we have $q\&q\&q=q\&q$ since $q\&q\leq q\&q$, hence $q\&q$  is idempotent. Let $a=q\&q$, then it follows that $p\leq a$. Otherwise, if $a< p\leq q\leq \hat{a}$, then $a\in S$,  which means that both $p$ and $q$ are in the same interval $[a,\hat{a}]$ with some $a\in S$,  a contradiction. Thus, one obtains that $$p\&q=p=p\wedge q.$$
\end{itemize}

Therefore, in both cases, it is shown that  \[p\& q=\begin{cases}
		a,&p,q\in [a,\hat{a}] \text{ for some } a\in S,\\
		p\wedge q,&\text{otherwise},
		\end{cases} \]
  as desired. 

$(3)\Lra(1)$: Given $p,q,u\in [0,1]$, it holds that 
 \begin{align*}
((p\wedge u)\& q)\vee (p\&(q\wedge u))
&=((p\&q)\wedge(u\&q))\vee((p\&q)\wedge(p\&u))\\
&=(p\&q)\wedge((u\&q)\vee(p\&u))\\
&=(p\&q)\wedge(u\&(p\vee q)).
\end{align*}
We check the equality 
$$(p\&q)\wedge (u\&(p\vee q))=(p\&q)\wedge u$$
in three cases:
\begin{enumerate}
 \item[(i)] If $p\vee q\leq u$, then both $p\leq u$ and $q\leq u$ hold.  We have that $p\&q\leq u\&(p\vee q)\leq u$. Thus, it follows that
 $$(p\&q)\wedge (u\&(p\vee q))=p\&q=(p\&q)\wedge u.$$
\item[(ii)] If $u<p\vee q$ and there is some $i\in I$ such that $a_i\leq u<p\vee q\leq b_i$, then we have that $p\&q\leq b_i\&b_i=a_i$ and $u\&(p\vee q)=a_i$. Thus, it follows that 
$$(p\&q)\wedge(u\&(p\vee q))=p\&q=(p\&q)\wedge u.$$
\item[(iii)] If $u<p\vee q$ but no $[a_i,b_i]$ contains $u$ and $p\vee q$ simultaneously, then it holds that 
$$(p\&q)\wedge(u\&(p\vee q))=(p\&q)\wedge(u\wedge(p\vee q))=(p\&q)\wedge u.$$
\end{enumerate}
Therefore, we have checked that the equality holds as desired. 
\end{proof}


\begin{prop}\label{approx}
Let $\&$ be a left continuous t-norm satisfying condition $(3)$ in Proposition \ref{E1}. We have  
 $$1=\sup\{p<1\mid p\&p=p\}.$$ 
\end{prop}
\begin{proof}
    Denote this family of pairwise disjoint closed intervals by $\{[a_i,b_i]\subseteq[0,1)\mid i\in I\}$. 
    Let $b=\sup_{i\in I} b_i$. If $b<1$, then each  $p\in(b,1)$ is idempotent since $p\&p=p\wedge p$; If $b=1$, since each $a_i$ is idempotent and $b_{i-1}<a_i<b_i<1$, then $\sup_{i\in I}a_i=1$.
\end{proof}

A $[0,1]$-category $(X,r)$ is {\em finite} if its underlying set $X$ is finite. All finite $[0,1]$-categories form a full subcategory of $\RCat$, denoted by $$\RfCat.$$ It is closed under finite products.

\begin{lem}
Let $(X,r)$ be a finite $[0,1]$-category. If $(X,r)\times{(-)}:\RfCat\lra{\RfCat}$ has a right adjoint, that is, the power $(Y,s)^{(X,r)}$ exists for every finite $[0,1]$-category $(Y,s)$, then it is given by  $(Y,s)^{(X,r)}=([(X,r),(Y,s)],d)$ where $d$ is the hom-functor from equation (\ref{exp hom}). 
\end{lem}
\begin{proof}Let $(X,r), (Y,s)$ be finite $[0,1]$-categories.
    Firstly, the power $(Y,s)^{(X,r)}$ has underlying set $[(X,r),(Y,s)]$ since \[\RfCat((X,r)\times{\mathbb{T},(Y,s)})\cong\RfCat(\mathbb{T},(Y,s)^{(X,r)}).\]
    So the power  $(Y,s)^{(X,r)}$ is the $[0,1]$-category $([(X,r),(Y,s)],d')$ for a suitable hom-functor $d'$.

    Secondly, on one hand, consider the counit on $(Y,s)$:
    \[\ev:(X,r)\times{([(X,r),(Y,s)],d')}\lra{(Y,s)},\]
    it follows that for every $x,y\in X$ and $f,g\in[(X,r),(Y,s)],$  $r(x,y)\wedge d'(f,g)\leq s(f(x),g(y))$, that is, $d'(f,g)\leq d(f,g)$ for all  $f,g\in[(X,r),(Y,s)].$
    One the other hand, let $f,g\in[(X,r),(Y,s)]$ and $p=d(f,g)$. We may assume that $f\neq g$, then consider a finite $[0,1]$-category $(Z,t)$ with underlying set $\{f,g\}$ and  $t(f,f)=t(g,g)=1, t(f,g)=p, t(g,f)=0$. Since $r(x,y)\wedge p\leq s(f(x),g(y))$ for all $x,y\in X$, the map $$h:(X,r)\times{(Z,t)}\lra{(Y,s)},\quad h(x,f)=f(x), h(x,g)=g(x)$$ is indeed a $[0,1]$-functor, then its transpose
    \[\hat{h}: (Z,t)\lra{([(X,r),(Y,s)],d')}\] is also a $[0,1]$-functor, hence $d(f,g)=p=t(f,g)\leq d'(\hat{h}(f),\hat{h}(g))=d'(f,g).$
    
    Thus, \[d'(f,g)=d(f,g)=\bv\{q\in Q\mid q\wedge r(x,y)\leq s(f(x),g(y))\text{ for all } x,y\in X\}\] for all $f,g\in[(X,r),(Y,s)].$
\end{proof}

\begin{thm}\label{RCat ccc}
Let $\&$ be a left continuous t-norm on $[0,1]$. Then the following statements are equivalent:
\begin{enumerate}
\item[{\rm (1)}] The category $\RCat$ is cartesian closed;
\item[{\rm (2)}] The category $\RfCat$ is cartesian closed;
\item[{\rm (3)}] Any one of the  three equivalent conditions in Proposition \ref{E1} holds. 
\end{enumerate}
\end{thm}
\begin{proof}
(1)$\Lra$(2): This follows from the facts that the terminal object $\mathbb{T}$, the products $(X,r)\times(Y,s)$ and the power objects $(Y,s)^{(X,r)}$ are all finite for all finite $[0,1]$-categories $(X,r)$ and $(Y,s)$. 

(2)$\Lra$(3): Let $X=\{x,y\}$ with $x\neq y$. For a given element $u\in[0,1]$,  let $r(x,x)=r(y,y)=1$, $r(x,y)=u$ and $r(y,x)=0$. Then $(X,r)$ is a $[0,1]$-category. 
    
For given $p,q\in[0,1]$, consider the following three maps from $X$ to $[0,1]$:
$$f(z)=r(x,z),\quad 
  g(z)=p\wedge r(x,z),\quad 
  h(z)=(p\&(q\wedge r(x,z)))\vee ((p\wedge u)\&(q\wedge r(y,z))).
$$
For all $x,y\in[0,1]$, let $$d_L(x,y)=\sup\{z\in[0,1]\mid x\& z\leq y\},$$ then $([0,1],d_L)$ is a $[0,1]$-category. Moreover, the maps $f$, $g$ and $h$ are all $[0,1]$-functors from $(X,r)$ to $([0,1],d_L)$. 

Let $Y=f(X)\cup g(X)\cup h(X)$ and equip it with the substructure of $([0,1],d_L)$, then we obtain a  $[0,1]$-category $(Y,s)$ with finite underlying set.  
It is easy to check that $f,g$ and $h$ are indeed $[0,1]$-functors from $(X,r)$ to $(Y,s)$.
	 
Moreover, consider the power object $(Y,s)^{(X,r)}=([(X,r),(Y,s)],d)$ in $\RfCat$. Then it holds that $p\leq d(f,g)$ since for all $a,b\in X$,
 \begin{align*}
 r(x,a)\& (p\wedge r(a,b))\leq p\wedge r(x,b)
        &\Longrightarrow p\wedge r(a,b)\leq d_L(f(a),g(b))=s(f(a),g(b)).
\end{align*}
  Similarly,  we can check that $q\leq d(g,h)$. 

Since the power $(Y,s)^{(X,r)}$ is a $[0,1]$-category, we have\begin{align*}
(p\&q)\wedge u&=(p\& q)\wedge r(x,y)\\
&\leq (d(f,g)\& d(g,h))\wedge r(x,y)\\
&\leq d(f,h)\wedge r(x,y)\\
&\leq d_L(f(x),h(y))\\
&=h(y).
\end{align*}
Notice that $h(y)=(p\&(q\wedge u))\vee ((p\wedge u)\&q)$, 
it follows that 
$$(p\&q)\wedge u\leq (p\&(q\wedge u))\vee ((p\wedge u)\&q).$$
Moreover, we always have that $p\&(q\wedge u)\leq p\&q$,  $p\&(q\wedge u)\leq u$, $(p\wedge u)\&q\leq p\&q$ and $(p\wedge u)\&q\leq u$. Thus we can see that 
$$(p\&(q\wedge u))\vee ((p\wedge u)\&q))\leq (p\&q)\wedge u.$$
Therefore, we have shown that 
$$(p\&q)\wedge u=((p\wedge u)\&q)\vee (p\&(q\wedge u))$$ for all $p,q,u$ in $[0,1]$, that is the first condition in Proposition \ref{E1}. 

(3)$\Lra$(1): By Proposition \ref{ccc con} straightforwardly. 
\end{proof}

\section{Cauchy complete real-enriched categories}
 In 1973 \cite{Lawvere1973}, Lawvere proposed that metric spaces could be seen as categories enriched over the quantale $([0,\infty]^{\rm op},+,0)$. Within this framework, Lawvere introduced a general concept of Cauchy completeness for enriched categories (via adjoint distributors). He then proved a remarkable theorem: for metric spaces viewed as enriched categories, this abstract categorical notion of Cauchy completeness is equivalent to the traditional, analytic completeness defined via Cauchy sequences. Further research of Cauchy completeness involves the categories enriched over a general quantale or even a quantaloid (see e.g., \cite{Hofmann2013a,Stubbe2005}).

 For real-enriched categories, Cauchy completeness is equivalent to the bicompleteness of Cauchy nets \cite{Zhang2024}.
 
 An element   $a$ in $X$ is a \emph{bilimit} of a net $\xlam$ in a $[0,1]$-category $(X,r)$ if for all $x\in X,$ $$r(a,x)=\sup_{\lam\in D}\inf_{\lam\leq\mu}r(x_{\mu},x)\text{ and }r(x,a)=\sup_{\lam\in D}\inf_{\lam\leq\mu}r(x,x_{\mu}).$$
It is clear that a net has at most one  bilimit up to isomorphism.
Some equivalent characterizations of bilimits of nets are given in \cite{Zhang2024}, even though in which real-enriched categories are based on a continuous t-norm, the following cited results are still valid for left-continuous t-norms.
\begin{lem}\label{bilim}{\rm (\cite{Zhang2024})}
    Let $\xlam$ be a  net in a $[0,1]$-category $(X,r)$ and $a$ be an element in $X$. Then the following statements are equivalent:
    \begin{enumerate}
    \item[{\rm (1)}] $a$ is a bilimit of $\xlam$.
    \item[{\rm (2)}] $\sup_{\lam\in D}\inf_{\lam\leq\mu}r(a,x_{\mu})=\sup_{\lam\in D}\inf_{\lam\leq\mu}r(x_{\mu},a)=1.$
    \item[{\rm (3)}] for all $\vep<1$, there is some $\lam\in D$ such that for all $\mu\geq\lam$, $\vep<r(a,x_\mu)$ and $\vep<r(x_\mu,a)$. 
    \end{enumerate} 
\end{lem}

 A net $\xlam$ in a $[0,1]$-category $(X,r)$  is \emph{Cauchy} if $$\sup_{\lam\in D}\inf_{\lam\leq\mu,\ga}r(x_{\mu},x_{\ga})=1,$$
or equivalently, $$\forall\vep<1, \exists\lam\in D,  \forall\mu,\nu\geq\lam, \vep<r(x_\mu,x_\nu).$$

A $[0,1]$-category $(X,r)$ is \emph{Cauchy complete} if every Cauchy net has a  bilimit.
  All  Cauchy complete $[0,1]$-categories and $[0,1]$-functors constitute a full subcategory of $\RCat$ that we denote as \[\CCom.\]

\begin{exmp}
    \begin{enumerate}
        \item Every finite $[0,1]$-category is trivially Cauchy complete, i.e. $\RfCat$ is a full sub category of $\CCom$.
        \item Let $(X,\leq)$ be a preordered set. The associated $[0,1]$-category $(X,r_\leq)$  is given by  
\[r_\leq(x,y)=\begin{cases}
		0,&x\nleq y,\\
		1,&x\leq y.
		\end{cases} \]
 A net $\xlam$ in $(X,r_\leq)$ is Cauchy if and only if it is eventually valued in isomorphic elements, that is, there is some $\lam\in D$ such that $x_{\mu}\approx x_{\nu}$ for all $\lam\leq\mu,\nu$. Thus, $(X,r_\leq)$ is always Cauchy complete.
        \item The $[0,1]$-category $([0,1],d_L)$  is Cauchy complete \cite{Wagner1997}.
    \end{enumerate}
\end{exmp}

\begin{lem}{\rm (\cite{Zhang2024})} \label{bilim pres} {\rm (1)} If a net $\xlam$ in a $[0,1]$-category $(X,r)$ has a bilimit, then $\xlam$ is a Cauchy net.

 {\rm (2)}   Each $[0,1]$-functor $f: (X,r)\lra (Y,s)$  preserves bilimits of Cauchy nets.
\end{lem}

\begin{prop}\label{Caucom prod}
   The category $\CCom$ has finite products. 
\end{prop}
\begin{proof}
    Let $(X,r)$ and $(Y,s)$ be Cauchy complete $[0,1]$-categories, it suffices to show that the product $(X\times Y,r\times{s})$ of $(X,r)$ and $(Y,s)$ is Cauchy complete.

    Suppose $\{(x_{\lam},y_{\lam})\}_{\lam\in D}$ is a Cauchy net in $(X\times Y,r\times{s})$, then $\xlam$ and $\ylam$ are Cauchy in $(X,r)$ and $(Y,s)$ respectively. Let $a\in X$ and $b\in Y$ denote the bilimit of $\xlam$ and $\ylam$. We claim that $(a,b)$ is the bilimit of $\{(x_{\lam},y_{\lam})\}_{\lam\in D}$ in $(X\times Y,r\times{s})$.

    Since $\sup_{\lam\in D}\inf_{\lam\leq\mu}r(a,x_{\mu})=\sup_{\lam\in D}\inf_{\lam\leq\mu}s(b,y_{\mu})=1$, it follows that \begin{align*}
        \sup_{\lam\in D}\inf_{\lam\leq\mu}r\times{s}((a,b),(x_{\mu},y_{\mu}))&=\sup_{\lam\in D}\inf_{\lam\leq\mu}(r(a,x_{\mu})\wedge s(b,y_{\mu}))\\
        &=\Big(\sup_{\lam\in D}\inf_{\lam\leq\mu}r(a,x_{\mu})\Big)\wedge\Big(\sup_{\lam\in D}\inf_{\lam\leq\mu}s(b,y_{\mu})\Big)\\
        &= 1.
    \end{align*}
    Similarly, $\sup_{\lam\in D}\inf_{\lam\leq\mu}r\times{s}((x_{\mu},y_{\mu}),(a,b))=1$. 
    By Lemma \ref{bilim}, $(a,b)$ is the bilimit of $\{(x_{\lam},y_{\lam})\}_{\lam\in D}$.
\end{proof}


\begin{thm}
    Let $\&$ be a left continuous t-norm. Then  
the category $\CCom$ is cartesian closed if and only if $\&$ satisfies any one of the three equivalent conditions in Proposition \ref{E1}.
\end{thm}
\begin{proof}
{\em Necessity}: Notice that the category $\RfCat$ is a full subcategory of $\CCom$ and is closed under finite products and power objects. Hence, the category $\RfCat$ is cartesian closed. Therefore, the necessity follows from Theorem \ref{ccc con}. 

{\em Sufficiency}: 
    Since the category $\RCat$ is cartesian closed, it suffices to show that the power $(Y,s)^{(X,r)}=([(X,r),(Y,s)],d)$ in $\RCat$ is Cauchy complete for all Cauchy complete $[0,1]$-categories $(X,r)$ and $(Y,s)$. 
    
    Let $\flam$ be a Cauchy net in $(Y,s)^{(X,r)}$, then for any given $p<1$, there is some $\lam\in D$ such that $$p<d(f_\mu,f_\nu)$$ for all $\mu\geq\lam$ and $\nu\geq\lam$. Thus, by equation (\ref{exp hom}), we have that
    $$p\wedge r(x,y)\leq s(f_\mu (x),f_\nu (y))$$ for all $x,y\in X$.  Particularly, let $x=y$, we obtain that
    $$p\leq s(f_\mu (x),f_\nu (x))$$ for all $x\in X$. Thus, $\{f_\lam (x)\}_{\lam\in D}$ is a Cauchy net in $(Y,s)$ and it has a bilimit, say $f(x)$. We claim that $$f:(X,r)\lra(Y,s),\quad x\mapsto f(x)$$ is a $[0,1]$-functor. In fact, for all $x,y\in X$, 
    \begin{align*}
    s(f(x),f(y))&=\sup_{\lam\in D}\inf_{\lam\leq\mu}s(f_\mu(x),f(y))\\
    &=\sup_{\lam\in D}\inf_{\lam\leq\mu}\sup_{\lam'\in D}\inf_{\lam'\leq\mu'}s(f_\mu(x),f_{\mu'}(y))\\
    &\geq\sup_{\lam\in D}\inf_{\lam\leq\mu}\sup_{\lam'\in D}\inf_{\lam'\leq\mu'}d(f_\mu,f_{\mu'})\wedge r(x,y)\\ 
    &=\sup_{\lam\in D}\inf_{\lam\leq\mu,\nu}d(f_\mu,f_{\nu})\wedge r(x,y)\\
    &=r(x,y).
    \end{align*}
   So, $f$ is a $[0,1]$-functor as claimed.

    Furthermore, we claim that $f$ is a bilimit of the Cauchy net $\flam$ in the power $(Y,s)^{(X,r)}$. 
    By Proposition 
    \ref{approx}, the set  $E=\{p\in[0,1)\mid p\&p=p\}$ satisfies that $\sup E=1$. We check that for all $q\in E$, there is some $\lam_q\in D$ such that for all $\mu\geq\lam_q$, 
    $$q\leq d(f_\mu,f).$$

     Since $f(x)$ is a bilimit of $\flamx$ in $(Y,s)$ for each $x\in X$, we have $$\sup_{\lam\in D}\inf_{\lam\leq\mu}s(f_\lam(x),f(x))=1.$$ That is, for each $p\in E$, there is some $\lam_p\in D$ such that  $p< s(f_\lam(x),f(x))$ for all $\mu\geq\lam_p$. In this case, for every $x,y\in X,$ since $p$ is idempotent, 
    \begin{align*}
        p\wedge r(x,y)&=p\& r(x,y)\\
        &\leq s(f_{\mu}(x),f(x))\& s(f(x),f(y))\\
        &\leq s(f_{\mu}(x),f(y)),
    \end{align*}  
    hence $p\leq d(f_{\mu},f)$ for all $\mu\geq\lam_p$. By $\sup E=1,$ we obtain that $$\sup_{\lam\in D}\inf_{\lam\leq\mu}d(f_{\mu},f)=1.$$ Similarly, we have that  $$\sup_{\lam\in D}\inf_{\lam\leq\mu}d(f,f_{\mu})=1.$$ By Lemma \ref{bilim}, $f$ is the bilimit of $\flam$ in $(Y,s)^{(X,r)}$ as desired.
 \end{proof}

\section{Yoneda complete and Smyth complete real-enriched categories}
    A net $\xlam$ in a $[0,1]$-category $(X,r)$ is \emph{forward Cauchy} \cite{Wagner1994,Wagner1997} if $$\sup_{\lam\in D}\inf_{\lam\leq\mu\leq\ga}r(x_{\mu},x_{\ga})=1.$$
The difference with Cauchy nets lies thus in the ordering imposed on the indices of this supinf.

An element $a$ in $X$ is a \emph{Yoneda limit} of a forward Cauchy net $\xlam$ if $$r(a,x)=\sup_{\lam\in D}\inf_{\lam\leq\mu}r(x_{\mu},x)$$ for all $x\in X$. A {Yoneda limit} of $\xlam$ is  denoted as $a=\mathrm{lim} x_\lam.$

A $[0,1]$-category $(X,r)$ is \emph{Yoneda complete} if every forward Cauchy net in $(X,r)$ has a Yoneda limit. For a $[0,1]$-functor $f:(X,r)\lra (Y,s)$ and a forward Cauchy net $\xlam$ in $(X,r)$, $\{f(x_\lam)\}_{\lam\in D}$ is clearly a forward Cauchy net in $(Y,s)$, the $[0,1]$-functor $f$ is \emph{Yoneda continuous} if it preserves the Yoneda limit of all forward Cauchy nets in $(X,r)$, that is, $f(\mathrm{lim}x_\lam)=\mathrm{lim} f(x_\lam)$ for all forward Cauchy nets $\xlam$ in $(X,r)$.
All   Yoneda complete $[0,1]$-categories and Yoneda continuous $[0,1]$-functors constitute a category 
\[\YCom.\]
 
\begin{rem}
   It is known that if we restrict the values of Yoneda complete $[0,1]$-categories on the two-point set $\{0,1\}$, we obtain a full subcategory of $\YCom$, which is isomorphic to a cartesian closed category $\mathsf{DCO}$ consisting of directed complete ordered sets and Scott continuous maps.
\end{rem}

Clearly, a Cauchy net $\xlam$ in a $[0,1]$-category $(X,r)$ is a forward Cauchy net, and an element $a\in X$ is a bilimit of $\xlam$ if and only if it is a Yoneda limit of $\xlam$. Thus, a Yoneda complete $[0,1]$-category is always Cauchy complete and  $\YCom$ is a (non-full) subcategory of $\CCom$.

 Smyth completeness was introduced by Smyth \cite{Smyth1988,Smyth1994} on quasi-uniform spaces, and it was extended to real-enriched categories \cite{Yu2023}.
 
A $[0,1]$-category is \emph{Smyth complete} if every forward Cauchy net has a bilimit. Since nets having bilimits must be Cauchy nets, a $[0,1]$-category is Smyth complete if and only if it is Cauchy complete and all forward Cauchy nets in it are Cauchy nets.
 
All  Smyth complete $[0,1]$-categories and $[0,1]$-functors constitute a category 
\[\SCom,\]
which is full in both  $\CCom$ and $\YCom$. Now, we obtain a chain of  categories: $$\RfCat\subseteq\SCom\subseteq\YCom\subseteq\CCom\subseteq\RCat.$$


\begin{prop}\label{product}
Both $\YCom$ and $\SCom$  are closed under finite products.
\end{prop}
\begin{proof}
    The proof is similar to that in Proposition \ref{Caucom prod}.
\end{proof}

 $\YCom$ is a concrete category over $\mathsf{Set}$,  and  the  terminal
object $\mathbb{T}$ is discrete. Thus, when   $\YCom$  is cartesian closed,  it also has function spaces. That is, for all real-enriched categories $(X, r)$ and $(Y, s)$ in $\YCom$, we can choose their power with $[(X,r)\ra(Y,s)]=\YCom((X,r),(Y,s))$  being underlying set and the evaluation morphism
$$\mathrm{ev}:X\times [(X,r)\ra(Y,s)]\lra Y,\quad \mathrm{ev}(x,f)=f(x).$$

 \begin{prop}\label{exp in yocq}
 Let  $\&$ be a left continuous t-norm which satisfies any one of the conditions in Proposition 
    \ref{E1}. Suppose $(X,r)$ and $(Y,s)$ are  $[0,1]$-categories, then \begin{itemize} 
 \item[{\rm (1)}]the $[0,1]$-category $([(X,r)\ra(Y,s)],d)$ is also Yoneda complete if both $(X,r)$ and $(Y,s)$ are;

\item[{\rm (2)}] the $[0,1]$-category $([(X,r),(Y,s)],d)$ is also Smyth complete if both $(X,r)$ and $(Y,s)$ are.
\end{itemize}
 \end{prop}
 
\begin{proof}
 (1)   The main technique is quite similar to that in \cite{Lai2025}, we present the sketch of the proof here for the convenience of the reader. In this case, since $\RCat$ is cartesian closed, the pair $([(X,r)\ra(Y,s)],d)$ is a subobject of the power of Yoneda complete $[0,1]$-categories $(X,r)$ and $(Y,s)$ in $\RCat$, hence it is indeed a $[0,1]$-category, then it suffices to show that it is Yoneda complete. Suppose  $\flam$ is a forward Cauchy net in $([(X,r)\ra(Y,s)],d)$. 
    
    Firstly, for every $x\in X$, $\flamx$ is a forward Cauchy net in $(Y,s)$ hence it has a Yoneda limit denoted as $f(x)$, one obtains a function $f:X\lra{Y}, x\mapsto f(x)=\lim f_{\lam}(x)$, it is also a Yoneda continuous $[0,1]$-functor from $(X,r)$ to $(Y,s)$ \cite[Theorem 4.2]{Lai2007}.

Secondly, for each idempotent element $p<1$, $p\& x=p\wedge x$ for every  $x\in [0,1]$ since the t-norm $\&$ satisfies condition (3) in Proposition 
    \ref{E1}. As shown in \cite[Lemma 5.8]{Lai2025}, for each $x\in X$, since the nets $\flam$ and $\flamx$ are eventually $p$-monotone in that sense, there is a $\lam_{0}\in D$ such that for all $\mu\geq\lam\geq\lam_0$, 
$$p\wedge d(f_\mu,g)\leq p\wedge d(f_\lam,g)$$  and
$$p\wedge s(f_\mu(x),g(y))\leq p\wedge s(f_\lam(x),g(y))$$ for all  $g\in[(X,r)\ra(Y,s)]$ and $x,y\in X$.
And as calculated similarly in \cite[Theorem 5.11]{Lai2025}, we have
    \[p\wedge d(f,g)=p\wedge\sup_{\lam\in D}\inf_{\lam\leq\mu} d(f_\mu,g)\] for all $g\in [(X,r)\ra (Y,s)]$.
    
Finally, since $\&$ is left continuous, and the set $E=\{p\text{ is idempotent }\mid p< 1 \}$ has  supremum $1$. It follows that \[d(f,g)=\sup_{p\in E}(p\wedge d(f,g))=\sup_{p\in E} (p\wedge  \sup_{\lam\in D}\inf_{\lam\leq\mu}d(f_\mu,g))=\sup_{\lam\in D}\inf_{\lam\leq\mu}d(f_\mu,g).\] for all $g\in [(X,r)\ra (Y,s)]$.

Hence $f$ is a Yoneda limit of $\flam$ in $([(X,r)\ra(Y,s)],d)$.

 (2)  For Smyth complete $[0,1]$-categories $(X,r),(Y,s)$. Notice that every $[0,1]$-functor $f:(X,r)\lra (Y,s)$ is Yoneda continuous by Lemma \ref{bilim pres}, hence $[(X,r),(Y,s)]=[(X,r)\ra(Y,s)]$. Suppose $\flam$ is a forward Cauchy net in $([(X,r),(Y,s)],d)$, denote the bilimit of $\flamx$ in $(Y,s)$ as $f(x)$ for each $x\in X$, one obtains a pointwise $[0,1]$-functor $f\in [(X,r),(Y,s)]$. Similarly, one has that  \[d(f,g)=\sup_{\lam\in D}\inf_{\lam\leq\mu}d(f_\mu,g)\] and  
    \[d(g,f)=\sup_{\lam\in D}\inf_{\lam\leq\mu}d(g,f_\mu)\] for all $g\in [(X,r),(Y,s)]$. That is, $f$ is the bilimit of $\flam$ in $([(X,r),(Y,s)],d)$, hence $([(X,r),(Y,s)],d)$ is a Smyth complete $[0,1]$-category.
\end{proof}

\begin{prop}\label{separate}
     Let $(X,r),(Y,s),(Z,t)$ be Yoneda complete $[0,1]$-categories. A $[0,1]$-functor $f:(X,r)\times (Y,s)\lra (Z,t)$ is Yoneda continuous if and only if it is Yoneda continuous separately.
\end{prop}
\begin{proof}
    The proof is similar to that of \cite[Proposition 3.7]{Lai2016}.
\end{proof}

\begin{prop}\label{evmap}  Let  $\&$ be a left continuous t-norm which satisfies any one of the conditions in Proposition 
    \ref{E1}. 
 Suppose $(X,r)$ and $(Y,s)$ are  Yoneda complete $[0,1]$-categories, then the evaluation map $\ev:(X,r)\times ([(X,r)\ra(Y,s)],d)\lra{(Y,s)}$ is a Yoneda continuous $[0,1]$-functor.
\end{prop}
\begin{proof}
    It is clear that $\ev$ is a $[0,1]$-functor and $\ev(-,f):(X,r)\lra{(Y,s)}$  is Yoneda continuous for each $f\in[(X,r)\ra(Y,s)]$. Given an $x\in X$  and a forward Cauchy net $\flam$ in $([(X,r)\ra(Y,s)],d)$, by Theorem \ref{exp in yocq}, the function $f:X\lra{Y}$ given by $f(x)=\lim f_{\lam}(x)$ is a Yoneda limit of $\flam$.  It follows that $\ev(x,f)=f(x)=\lim f_{\lam}(x)$, hence $\ev(x,-)$ is Yoneda continuous. By Proposition \ref{separate}, $\ev$ is Yoneda continuous.
\end{proof}

\begin{prop}\label{transp}   Let  $\&$ be a left continuous t-norm which satisfies any one of the conditions in Proposition 
    \ref{E1}.  Suppose $(X,r)$, $(Y,s)$ and $(Z,t)$ are all Yoneda complete $[0,1]$-categories and $f:(X,r)\times{(Z,t)}\lra{(Y,s)}$ is a Yoneda continuous $[0,1]$-functor, then $\hat{f}:(Z,t)\lra{([(X,r)\ra(Y,s)],d)}$ is also a Yoneda continuous $[0,1]$-functor, where $\hat{f}(z)=f(-,z)$ for each $z\in Z$.
\end{prop}
\begin{proof}
    It is easily verified that $\hat{f}$ is a $[0,1]$-functor. Given a forward Cauchy net $\zlam$ in $(Z,t)$ whose Yoneda limit denoted by $a\in Z$, then $\{f(-,z_\lam)\}_{\lam\in D}$ is a forward Cauchy net in $([(X,r)\ra(Y,s)],d)$ with a pointwise   Yoneda limit $g$, which is given by 
    $$g(x)=\lim f(x,z_\lam).$$ 
 Since $f$ is Yoneda continuous, it follows that $\hat{f}(a)(x)=f(x,a)=\lim f(x,z_\lam)=g(x)$ for all $x\in X$, that is,  $\hat{f}(a)=g$. The $[0,1]$-functor $\hat{f}$ is Yoneda continuous as desired.
\end{proof}

Combining the results of Proposition \ref{product}, Proposition \ref{exp in yocq}, Proposition \ref{evmap} and Proposition \ref{transp}, we have the following:
\begin{prop}\label{ycoq is ccc}
     Let  $\&$ be a left continuous t-norm which satisfies any one of the conditions in Proposition 
    \ref{E1}.  Then the categories $\YCom$ and $\SCom$ are cartesian closed.
\end{prop}


\begin{thm} Let $\&$ be a left continuous t-norm. The following statements are equivalent:  \begin{enumerate} 
 \item[{\rm (1)}] The category $\SCom$ is cartesian closed;
 \item[{\rm (2)}] The category $\YCom$ is cartesian closed.
 \item[{\rm (3)}] $\&$ satisfies any one of the three equivalent conditions in Proposition \ref{E1}.
\end{enumerate} \end{thm}

\begin{proof}  It suffices to show that $\RfCat$ is cartesian closed, this follows from  finite $[0,1]$-categories are trivially Yoneda complete and Smyth complete, and $\RfCat$ is closed under finite products and power objects.
\end{proof}

 
 \section{Conclusion}
Real-enriched categories are often treated as quantitative ordered sets in Lawvere's sense. It is known that if $\&$ is a continuous t-norm on $[0,1]$, the category $\qcat$ is catesian closed if and only if the quantale $([0,1],\&,1)$ is locally cartesian closed, that is, $\&=\wedge$. Dropping the continuity of the binary operation $\&$, we describe all left continuous operation $\&$ such that $\qcat$ is cartesian closed. In fact, we show that all statements in the below are equivalent:   
\begin{enumerate}
\item[{\rm (1)}] The left continuous t-norm $\&$ satisfies any one of the three equivalent conditions in Proposition \ref{E1};
\item[{\rm (2)}] $\RCat$, consisting of all $[0,1]$-categories and $[0,1]$-functors, is cartesian closed;
\item[{\rm (3)}] $\RfCat$, consisting of all finite $[0,1]$-categories and $[0,1]$-functors, is cartesian closed;
\item[{\rm (4)}] $\CCom$, consisting of all Cauchy complete $[0,1]$-categories and $[0,1]$-functors, is cartesian closed;
\item[{\rm (5)}] $\YCom$, consisting of all Yoneda complete $[0,1]$-categories and Yoneda continuous $[0,1]$-functors, is cartesian closed;
\item[{\rm (6)}] $\SCom$, consisting of all Smyth complete $[0,1]$-categories and $[0,1]$-functors, is cartesian closed.
\end{enumerate}

\bibliography{qingzhuRef}

@Book{Adamek1990,
  title     = {Abstract and Concrete Categories: The Joy of Cats},
  publisher = {Wiley},
  year      = {1990},
  author    = {Ad{\'a}mek, Ji{\v r}{\'i} and Herrlich, Horst and Strecker, George E.},
  address   = {New York},
  isbn      = {9780471609223},
  lccn      = {89014835},
  url       = {http://www.tac.mta.ca/tac/reprints/articles/17/tr17abs.html},
}

@Article{Clementino2009,
  author  = {Clementino, Maria Manuel and Hofmann, Dirk},
  title   = {Lawvere Completeness in Topology},
  journal = {Applied Categorical Structures},
  year    = {2009},
  volume  = {17},
  number  = {2},
  pages   = {175--210},
  issn    = {1572-9095},
  doi     = {10.1007/s10485-008-9152-5},
  url     = {http://dx.doi.org/10.1007/s10485-008-9152-5},
}

@Article{Clementino2009b,
  author  = {Clementino, Maria Manuel
and Hofmann, Dirk and Stubbe, Isar},
  title   = {Exponentiable Functors Between Quantaloid-Enriched Categories},
  journal = {Applied Categorical Structures},
  year    = {2009},
  volume  = {17},
  number  = {1},
  pages   = {91--101},
  issn    = {1572-9095},
  doi     = {10.1007/s10485-007-9104-5},
  url     = {https://doi.org/10.1007/s10485-007-9104-5},
}

@article{Hofmann2012,
title = {A duality of quantale-enriched categories},
journal = {Journal of Pure and Applied Algebra},
volume = {216},
number = {8},
pages = {1866-1878},
year = {2012},
note = {Special Issue devoted to the International Conference in Category Theory `CT2010'},
issn = {0022-4049},
doi = {https://doi.org/10.1016/j.jpaa.2012.02.024},
url = {https://www.sciencedirect.com/science/article/pii/S0022404912000631},
author = {Dirk Hofmann and Paweł Waszkiewicz}
}

@article{Flagg1997,
title = {Continuity spaces: Reconciling domains and metric spaces},
journal = {Theoretical Computer Science},
volume = {177},
number = {1},
pages = {111-138},
year = {1997},
issn = {0304-3975},
doi = {https://doi.org/10.1016/S0304-3975(97)00236-3},
url = {https://www.sciencedirect.com/science/article/pii/S0304397597002363},
author = {Bob Flagg and Ralph Kopperman}
}

@Article{Hofmann2011,
  author   = {Hofmann, Dirk and Stubbe, Isar},
  journal  = {Topology and its Applications},
  title    = {Towards {Stone} duality for topological theories},
  year     = {2011},
  issn     = {0166-8641},
  number   = {7},
  pages    = {913--925},
  volume   = {158},
  doi      = {10.1016/j.topol.2011.01.010},
  url      = {http://www.sciencedirect.com/science/article/pii/S0166864111000113},
}

@Book{Klement2000,
  title     = {Triangular Norms},
  publisher = {Springer},
  year      = {2000},
  author    = {Klement, Erich Peter and Mesiar, Radko and Pap, Endre},
  volume    = {8},
  series    = {Trends in Logic},
  address   = {Dordrecht},
  isbn      = {978-90-481-5507-1},
  doi       = {10.1007/978-94-015-9540-7},
  url       = {http://link.springer.com/10.1007/978-94-015-9540-7},
}

@Article{Lai2007,
  author   = {Hongliang Lai and Dexue Zhang},
  title    = {Complete and directed complete {$\Omega$}-categories},
  journal  = {Theoretical Computer Science},
  year     = {2007},
  volume   = {388},
  pages    = {1--25},
  issn     = {0304-3975},
  doi      = {10.1016/j.tcs.2007.09.012},
  url      = {http://www.sciencedirect.com/science/article/pii/S0304397507006792},
}

@Article{Lawvere1973,
  author  = {Lawvere, Francis William},
  title   = {Metric spaces, generalized logic and closed categories},
  journal = {Rendiconti del Seminario Mat\'{e}matico e Fisico di Milano},
  year    = {1973},
  volume  = {43},
  pages   = {135--166},
  url     = {http://tac.mta.ca/tac/reprints/articles/1/tr1abs.html},
}

@Article{Stubbe2005,
  author  = {Stubbe, Isar},
  journal = {Theory and Applications of Categories},
  title   = {Categorical structures enriched in a quantaloid: categories, distributors and functors},
  year    = {2005},
  number  = {1},
  pages   = {1--45},
  volume  = {14},
  url     = {http://www.tac.mta.ca/tac/volumes/14/1/14-01abs.html},
}

@PhdThesis{Wagner1994,
  author  = {Wagner, Kim Ritter},
  title   = {Solving Recursive Domain Equations with Enriched Categories},
  school  = {Carnegie Mellon University},
  year    = {1994},
  address = {Pittsburgh},
}

@Article{Wagner1997,
  author  = {Wagner, Kim Ritter},
  title   = {Liminf convergence in {$\Omega$}-categories},
  journal = {Theoretical Computer Science},
  year    = {1997},
  volume  = {184},
  number  = {1-2},
  pages   = {61--104},
  issn    = {0304-3975},
  doi     = {http://dx.doi.org/10.1016/S0304-3975(96)00223-X},
  url     = {http://www.sciencedirect.com/science/article/pii/S030439759600223X},
}

@Article{Bonsangue1998,
  author  = {Bonsangue, M. M. and van Breugel, F. and Rutten, J. J. M. M.},
  title   = {Generalized metric spaces: Completion, topology, and powerdomains via the {Yoneda} embedding},
  journal = {Theoretical Computer Science},
  year    = {1998},
  volume  = {193},
  number  = {1},
  pages   = {1--51},
  issn    = {0304-3975},
  doi     = {http://dx.doi.org/10.1016/S0304-3975(97)00042-X},
  url     = {http://www.sciencedirect.com/science/article/pii/S030439759700042X},
}

@Article{Hofmann2013a,
  author   = {Hofmann, Dirk and Reis, C. D.},
  title    = {Probabilistic metric spaces as enriched categories},
  journal  = {Fuzzy Sets and Systems},
  year     = {2013},
  volume   = {210},
  pages    = {1--21},
  issn     = {0165-0114},
  doi      = {http://dx.doi.org/10.1016/j.fss.2012.05.005},
  url      = {http://www.sciencedirect.com/science/article/pii/S0165011412002242},
}

@Article{Lai2017,
  author   = {Lai, Hongliang and Shen, Lili},
  title    = {Fixed points of adjoint functors enriched in a quantaloid},
  journal  = {Fuzzy Sets and Systems},
  year     = {2017},
  volume   = {321},
  pages    = {1--28},
  doi      = {http://dx.doi.org/10.1016/j.fss.2016.12.001},
  issn     = {0165-0114},
  url      = {http://www.sciencedirect.com/science/article/pii/S0165011416304031},
}

@Article{Lai2016,
  author   = {Lai, Hongliang and Zhang, Dexue},
  title    = {Closedness of the category of liminf complete fuzzy orders},
  journal  = {Fuzzy Sets and Systems},
  year     = {2016},
  volume   = {282},
  pages    = {86--98},
  issn     = {0165-0114},
  doi      = {https://doi.org/10.1016/j.fss.2014.10.031},
  url      = {http://www.sciencedirect.com/science/article/pii/S0165011414004874},
}

@Article{Zhang2024,
  author  = {Zhang, Dexue},
  journal = {arXiv:2403.09716},
  title   = {Introductory notes on real-enriched categories},
  year    = {2024},
  url     = {https://arxiv.org/abs/2403.09716},
}

@Article{Clementino2003b,
  author  = {Maria Manuel Clementino and Dirk Hofmann and Walter Tholen},
  journal = {Theory and Applications of Categories},
  title   = {Exponentiability in categories of lax algebras},
  year    = {2003},
  number  = {15},
  pages   = {337--352},
  volume  = {11},
  url     = {https://api.semanticscholar.org/CorpusID:18294255},
}

@Article{Clementino2006,
  author   = {Maria Manuel Clementino and Dirk Hofmann},
  journal  = {Topology and its Applications},
  title    = {Exponentiation in $\mathbf{V}$-categories},
  year     = {2006},
  issn     = {0166-8641},
  note     = {Special Issue: Aspects of Contemporary Topology},
  number   = {16},
  pages    = {3113-3128},
  volume   = {153},
  doi      = {https://doi.org/10.1016/j.topol.2005.01.038},
  url      = {https://www.sciencedirect.com/science/article/pii/S0166864105001902},
}

@Article{Lai2025,
  author   = {Hongliang Lai and Qingzhu Luo},
  journal  = {Fuzzy Sets and Systems},
  title    = {Cartesian closed and stable subconstructs of [0,1]-{Cat}},
  year     = {2025},
  issn     = {0165-0114},
  pages    = {109243},
  volume   = {503},
  doi      = {https://doi.org/10.1016/j.fss.2024.109243},
  url      = {https://www.sciencedirect.com/science/article/pii/S0165011424003890},
}

@Article{Stubbe2025,
  author  = {Isar Stubbe and Junche Yu},
  journal = {arXiv:2501.03942},
  title   = {When is {{\sf Cat}($\mathcal{Q}$)} cartesian closed?},
  year    = {2025},
  url     = {
https://doi.org/10.48550/arXiv.2501.03942},
}

@InProceedings{Smyth1988,
  author    = {Smyth, M. B.},
  booktitle = {Mathematical Foundations of Programming Language Semantics},
  title     = {Quasi-uniformities: Reconciling domains with metric spaces},
  year      = {1988},
  address   = {Berlin, Heidelberg},
  editor    = {Main, M. and Melton, A. and Mislove, M. and Schmidt, D.},
  pages     = {236--253},
  publisher = {Springer Berlin Heidelberg},
  isbn      = {978-3-540-38920-0},
}

@Article{Smyth1994,
  author  = {Smyth, M. B.},
  journal = {Journal of the London Mathematical Society},
  title   = {Completeness of Quasi-Uniform and Syntopological Spaces},
  year    = {1994},
  number  = {2},
  pages   = {385-400},
  volume  = {49},
  doi     = {https://doi.org/10.1112/jlms/49.2.385},
  eprint  = {https://londmathsoc.onlinelibrary.wiley.com/doi/pdf/10.1112/jlms/49.2.385},
  url     = {https://londmathsoc.onlinelibrary.wiley.com/doi/abs/10.1112/jlms/49.2.385},
}

@Article{Yu2023,
  author        = {Junche Yu and Dexue Zhang},
  title         = {Smyth complete real-enriched categories},
  year          = {2023},
  journal = {arXiv:2311.18191},
  url           = {https://arxiv.org/abs/2311.18191},
}
\end{document}